\documentclass[12pt]{aipproc}
\layoutstyle{6s}
\usepackage{amssymb}
\usepackage{amsmath}
\usepackage{mathrsfs}
\usepackage{amsthm}
\usepackage[all]{xy}
\newtheorem{theorem}{Theorem}

\def\RR{\mathbb{R}}

\def\ttheta{\Psi}
\def\H{\mathscr{H}}
\def\y{\xi}

\def\MM{\widetilde{M}}
\def\MMM{\widehat{M}}
\def\V{\mathfrak{V}}
\def\ppi{\tilde{\pi}}
\def\pppi{\hat{\pi}}

\def\ggamma{\tilde{\Gamma}}
\def\gggamma{\widehat{\Gamma}}
\def\xx{\tilde{x}}
\def\xxx{\hat{x}}
\def\nnabla{\widetilde{\nabla}}
\def\w{\mathbf{w}}
\begin{document}
\newtheorem*{proposition}{Proposition}
\theoremstyle{remark}
\newtheorem*{remarks}{Remarks}
\newtheorem*{remark}{Remark}
\newtheorem*{definition}{Definition}
\title{Projective Connections and the Algebra of Densities}
\author{Jacob George}{address = {School of Mathematics, University of Manchester, Oxford Road, Manchester M13 9PL, UK}, email = {jgeorge@maths.man.ac.uk}}, 
\begin{abstract}
Projective connections first appeared in Cartan's papers in the
1920's.  Since then they have resurfaced periodically in, for
example, integrable systems and perhaps most recently in the context
of so called projectively equivariant quantisation.  We recall the
notion of projective connection and describe its relation with the
algebra of densities on a manifold.  In particular, we construct a
Laplace-type operator on functions using a Thomas projective
connection and a symmetric contravariant tensor of rank 2 (`upper
metric').
\end{abstract}
\maketitle
\section{Introduction}
This paper is concerned with the geometry of differential operators
on a manifold and their relation with projective connections.  The
notion of projective connection is an old one, first appearing in
Cartan's papers of the 1920s and then in various modified forms
throughout the 20th century. They have surfaced periodically in
mathematical physics, in particular in integrable systems and more
recently in projectively equivariant quantisation. This paper
establishes that the algebra of densities introduced in \cite{oloii}
and projective connections on a given manifold are fundamentally
linked.  We recall the notion of a projective connection; in fact
there are two distinct but related notions which will be detailed
here. We define the manifold $\MMM$ (cf. $\hat{M}$ in \cite{oloii})
for which the algebra of densities $\V(M)$ may be interpreted as a
subalgebra of the algebra of smooth functions $C^{\infty}(\MMM)$.
Having defined Thomas' manifold $\MM$ (see \cite{thomas3}), we
explicitly construct a diffeomorphism $F : \MMM \to \MM$. By means
of this, we show that a projective connection on $M$ gives rise
canonically to a linear connection on $\MMM$.  Then we consider a
manifold equipped with a Thomas projective connection and a
symmetric contravariant tensor field of rank 2 (which may be viewed
as an `upper metric', though it need not be non-degenerate). From
these data, we construct an upper connection (contravariant
derivative) on the bundle of volume forms and an invariant
Laplace-type operator acting on functions, a `projective Laplacian'.
This leads to some results regarding differential operators and
brackets on the algebra $\V(M)$ when $M$ is equipped with a
projective connection.
\section{Projective connections: a recollection}
Here we present a natural definition of projective connection and
give the correspondence to the Thomas-Weyl-Veblen definition of
projective equivalence classes.  Most of this account may be
extracted, with a little work, from \cite{hermann}. Recall that a
map $\phi:\RR^{n}\to\RR^{n}$ is said to be \emph{projective} or
\emph{fractional linear} if it is of the form
\begin{equation*}
\phi(v) = \frac{\alpha(v) + \beta}{\gamma(v) + \delta},\; \mbox{
where }\alpha\in GL_{n}(\RR),\beta\in\RR^{n},\gamma\in(\RR^{n})^{*},
\delta\in(\RR), \det\left(\begin{array}{c|c}\alpha&\beta\\\hline
\gamma&\delta\end{array}\right)\neq 0
\end{equation*}
With this in mind, we make the following natural definition.
\begin{definition}\label{heurconn}
A \emph{projective connection} on a vector bundle $E\to M$ is an
Ehresmann connection such that the associated parallel transport
induces a projective map between fibres.
\end{definition}
This intuitive definition is a reformulation, in modern language, of
Cartan's idea as proposed in his seminal paper \cite{Cartan}.  In
the special case of the tangent bundle $TM \to M$, projective
connections are closely related to the following notion.
\begin{definition}\label{twconn}
Two torsion free linear connections on $TM \to M$ are said to be
\emph{projectively equivalent} if any of the following hold.
\begin{enumerate}
\item $\nabla\sim\bar{\nabla}$ if they define the same geodesics up to reparametrisation.\label{1}
\item $\nabla\approx\bar{\nabla}$ if there is a $1$-form $\vartheta$
s.t. $\bar{\nabla}_{X}Y - \nabla_{X}Y = \vartheta(X)Y +
\vartheta(Y)X$ for all vector fields $X,Y$.\label{2}
\item $\nabla\simeq\bar{\nabla}$ if their coefficients $\Gamma^{k}_{ij}$ and $\bar{\Gamma}^{k}_{ij}$ satisfy\label{3}
\begin{equation*}
\Pi^{k}_{ij}:=\Gamma^{k}_{ij} - \frac{1}{n+1}(\delta^{k}_{i}\Gamma^{l}_{lj} + \delta^{k}_{j}\Gamma^{l}_{il})= \bar{\Gamma}^{k}_{ij} - \frac{1}{n+1}(\delta^{k}_{i}\bar{\Gamma}^{l}_{lj} + \delta^{k}_{j}\bar{\Gamma}^{l}_{il}) =: \bar{\Pi}^{k}_{ij}
\end{equation*}
\end{enumerate}
\end{definition}
This would not be a well defined notion had these conditions not
been equivalent.  Details of the proof of
(\ref{1})$\Leftrightarrow$(\ref{2}) can be found in \cite{hermann};
a brute force calculation yields
(\ref{2})$\Leftrightarrow$(\ref{3}). Before relating these two
notions, we must rephrase the former analytically.

Recall that an Ehresmann connection in a fibre bundle $\pi:E\to M$
is a distribution $\H$ of horizontal linear subspaces of $TE$ such
that $TE = \H\oplus\ker\pi_{*}$.  By \emph{horizontal curves} we
will mean integral curves of $\H$. We express $\H$ as the kernel of
a collection of $1$-forms on $E$. In terms of local coordinates
$x^{i}$ on the base and $\y^{a}$ on the fibres, these forms may be
written (after a careful choice of coordinates) in the form
$\ttheta^{a} = d\y^{a} + f^{a}_{i}dx^{i}$ for some functions
$f^{a}_{i}$ on $M$. These functions together constitute a local
connection $1$-form on $M$. The condition that a curve $\sigma$ be
horizontal is then simply that $\ttheta^{a}(\sigma') = 0$ for each
$a$, or explicitly in terms of $f^{a}_{i}$,
\begin{equation}\label{annihilate}
\frac{d}{dt}\y^{a}(\sigma(t)) = f^{a}_{i}\frac{d}{dt}x^{i}(\sigma(t)).
\end{equation}
In this context, we require that the parallel translation defines a
projective map between fibres.  This is equivalent to requiring that
locally, the horizontal curves are projective flows on the typical
fibre, that is, of the form
\begin{equation*}
\sigma(t) = \frac{\alpha(\sigma) + \beta}{\gamma(\sigma) + \delta}
\end{equation*}
with the aforementioned non-degeneracy condition on the
coefficients. Differentiating this, we have the equation in fibre
coordinates
\begin{equation*}
\frac{d}{dt}\y^{a}(\sigma(t)) = A^{a}(t) + B^{a}_{b}(t)\y^{b}(\sigma(t)) + C_{b}\y^{a}(\sigma(t))\y^{b}(\sigma(t))
\end{equation*}
for some coefficient functions $A^{a}$, $B^{a}_{b}$ and $C_{b}$.
Comparing this with (\ref{annihilate}) we come to a description of
projective connection in a vector bundle $E\to M$ as an Ehresmann
connection defined by the annihilating $1$-forms
\begin{equation*}
\ttheta^{a} = d\y^{a} - (\phi^{a}_{i}(x) + \psi^{a}_{bi}(x)\y^{b} +
\eta_{bi}(x)\y^{a}\y^{b})dx^{i}
\end{equation*}
for some functions $\phi^{a}_{i}$, $\psi^{a}_{bi}$ and $\eta_{bi}$.
\begin{remark} With respect to (linear) transition functions in the vector bundle
$E$, the forms $\psi^{a}_{bi}dx^{i}$ define a connection $1$-form
for a linear connection in $E$, while $\psi^{a}_{bi}dx^{i}$ and
$\eta_{bi}dx^{i}$ define $1$-forms taking values in sections of $E$
and the dual bundle $E^{*}$ respectively. This will be important in
the
following section.  %In the case $E = T(M)$, $\phi_{i}^{j}dx^{i}$ and
%$\eta_{ji}dx^{i}$ define $1$-forms taking values in vectors and
%covectors respectively.
\end{remark}
In the case of the tangent bundle, there is a distinguished choice
of $\phi^{j}_{i}$, namely $\delta^{j}_{i}$.  This gives $\ttheta^{i}
= d\y^{i} - dx^{i} - \omega^{i}_{j}\y^{j} -
\omega^{0}_{j}\y^{i}\y^{j}$ where now $\y^{i}$ are the fibre
coordinates naturally related to the local coordinates $x^{i}$.
Projective connections also have a related notion of geodesics.
These are given in this case by the equations
\begin{equation*}
\frac{d^{2}\sigma^{k}}{dt^{2}} = \frac{d\sigma^{a}}{dt} +
\Gamma^{k}_{ij}\frac{d\sigma^{i}}{dt}\frac{d\sigma^{j}}{dt} +
\gamma^{0}_{ij}\frac{d\sigma^{i}}{dt}\frac{d\sigma^{j}}{dt}\frac{d\sigma^{k}}{dt}
\end{equation*}
where $\omega^{k}_{i} =: \Gamma^{k}_{ij}dx^{j}$ and $\omega^{0}_{i}
=: \gamma^{0}_{ij}dx^{j}$.
 Having established this notion, we may mimic our previous definition and analogously define \emph{projective equivalence classes of
projective connections}.

For arbitrary Ehresmann connections, there is a natural notion of curvature which in this case
amounts to a collection of 2-forms:
\begin{equation*}
\Omega^{i}_{j} = d\omega^{i}_{j} -
\omega^{i}_{k}\wedge\omega^{k}_{j} - \omega^{0}_{j}\wedge dx^{i} -
\delta^{i}_{j}\omega^{0}_{k}\wedge dx^{k}\mbox{ and } \Omega^{0}_{i}
= \omega^{0}_{j}\wedge\omega^{j}_{i}.
\end{equation*}
%These are exactly the failure of the distribution $\H$ to be integrable, or in other words, the failure of $d\ttheta^{i}$ to lie in the ideal generated by $\ttheta^{j}\wedge\ttheta^{k}$.
\begin{definition}
Define the numbers $A^{i}_{jkl}$ and $A_{jk}$ by $\Omega^{i}_{j} =:
A^{i}_{jkl}dx^{k}\wedge dx^{l}$ and $A_{jk}:= A^{i}_{jki}$.  A
projective connection is said to be \emph{normal} if $A_{jk} =0$.
\end{definition}
We also consider the Ricci tensor $R_{jk}$, where $d\omega^{i}_{j} -
\omega^{i}_{p}\wedge\omega^{p}_{j}$ equals $R^{i}_{jkl}dx^{k}\wedge
dx^{l}$ and $R_{jk}:=R^{i}_{jki}$.

The notion of normality was introduced by Cartan and derives its utility from the following
observations. Let $(\omega^{i}_{j}, \omega^{0}_{j})$ define a projective connection in $TM$.
\begin{itemize}
\item Given that the Ricci tensor $R_{ij}$ is symmetric, the projective connection is normal if and only if
$\omega^{0}_{j} = \frac{2}{n-1}R_{jk}dx^{k}$.
\item Given any linear connection $\omega^{i}_{j}$ in the tangent bundle, there is another with symmetric Ricci tensor with the same geodesics up to reparametrisation.
\item If $(\omega^{i}_{j}, \omega^{0}_{j})$, $(\bar{\omega}^{i}_{j},\bar{\omega}^{0}_{j})$ are two normal
projective connections in the tangent bundle whose associated linear
connections $\omega^{i}_{j}$ and $\bar{\omega}^{i}_{j}$ are
projectively equivalent, then $(\omega^{i}_{j}, \omega^{0}_{j})$ and
$(\bar{\omega}^{i}_{j},\bar{\omega}^{0}_{j})$ are projectively
equivalent.
\end{itemize}
Proofs of these facts may be be found in \cite{hermann}.  Given
these, we have a one-to-one correspondence between projective
equivalence classes of linear connections and projective equivalence
classes of projective connections on any given manifold $M$.  All of
the notions discussed throughout this paper also make sense on
supermanifolds.

Bearing these relations in mind, we will adopt from now on the
(classical) terminology and call a projective equivalence class of
linear connections on a manifold a \emph{Thomas projective
connection}.
%Through this proposition, we see the relation between Definition \ref{heurconn} and Definition
%\ref{twconn}:  To generate a projective equivalence class from a projective connection
%$(\omega^{a}_{b},\omega^{0}_{b})$, simply take $[(\omega^{a}_{b})]$. Conversely, given a projective
%equivalence class (and a representative) $[(\omega^{a}_{b})]$, there is a canonical normal
%projective connection, namely $\left(\omega^{a}_{b},\frac{2}{n-1}R_{jk}dx^{k}\right)$.
\section{Relations with densities}
\begin{definition}A \emph{density of weight $\lambda$} is a formal expression of the form $\mathbf{\phi}=
\phi(x)(Dx)^{\lambda}$ defined in local coordinates, $Dx$ being the
associated local volume form and $\lambda\in\RR$.  There is a
natural notion of multiplication for densities given by
$\mathbf{\phi}\cdot\mathbf{\chi} =
(\phi(x)(Dx)^{\lambda_{1}})\cdot(\chi(x)Dx^{\lambda_{2}}) =
\phi(x)\chi(x)(Dx)^{\lambda_{1} + \lambda_{2}}$.  The \emph{algebra
of densities} denoted $\V(M)$, is the algebra of finite formal sums
$\sum_{\lambda}\phi_{\lambda}(x)(Dx)^{\lambda}$ of densities with
the multiplication defined previously (see \cite{oloii}).
\end{definition}

A $\lambda$-density may equally be thought of as a `function'
$\phi(x)$ which under a change of coordinates, picks up the modulus
of the Jacobian of the transformation to the $\lambda$-th power as a
factor. Having defined densities, it is natural to define a linear
operator on $\V(M)$, namely the weight operator $\w$. This is
defined as having each $\lambda$-density as a $\lambda$-eigenvector,
that is $\w(\phi(x)(Dx)^{\lambda}) =\lambda\phi(x)(Dx)^\lambda$.

Let $\dim M = n$.  We now define two $(n+1)$-dimensional manifolds
$\MMM$ and $\MM$ which are fibre bundles over $M$ and significant
with respect to the algebra of densities and projective connections
respectively.
\begin{description}
\item[$\MMM:$ Densities as functions.] A density $\sum \phi_{\lambda}(Dx)^{\lambda}$ may be interpreted as a function on the manifold $\MMM$ defined
in some sense as the `strictly positive half' of the determinant
bundle (see below).  Denoting the fibre coordinate by $t(> 0)$, the
algebra of densities is the subalgebra of $C^{\infty}(\MMM)$
consisting of functions of the form $\sum
\phi_{\lambda}t^{\lambda}$.
\item[$\MM:$ Thomas projective connections as linear connections.] In \cite{thomas2}, T. Y. Thomas details the construction of a
manifold $\MM$ from a given manifold $M$ such that any projective
equivalence class $\Pi^{k}_{ij}$ on $M$ gives rise canonically to
linear connection coefficients $\ggamma^{k}_{ij}$ on $\MM$.  From
the viewpoint of the previous section, each projective equivalence
class of projective connections gives rise to a linear connection on
$\MM$.  If $x^{0},\dots,x^{n}$ are local coordinates in some
neighbourhood of $\MM$, these are given by
\begin{eqnarray}
\ggamma^{k}_{ij} = \ggamma^{k}_{ji} = \Pi^{k}_{ij} &&\mbox{for}\;
i,j,k = 1,\dots,n\label{t1}\\
\ggamma^{k}_{i0} = \ggamma^{k}_{0i} = -\frac{\delta^{k}_{i}}{n+1}
&&\mbox{for}\;i,k = 0,\dots,n\label{t2}\\
\ggamma^{0}_{ij} = \ggamma^{0}_{ji} = \frac{n+1}{n-1}\left(\frac{\partial\Pi^{r}_{ij}}{\partial
x^{r}} - \Pi^{r}_{si}\Pi^{s}_{rj}\right) &&\mbox{for}\; i,j = 1,\dots,n.\label{t3}
\end{eqnarray}
\end{description}
The manifolds $\MM$ and $\MMM$ are defined locally in terms of the charts on $M$.  Let us for a
moment denote local coordonates on $\MMM$ by $t, \xxx^{1},\dots,\xxx^{n}$ and on $\MM$ by
$\xx^{0},\dots,\xx^{n}$.  Then if $x'^{i} = f^{i}(x^{1},\dots,x^{n})$ is a coordinate
transformation on $M$, define
\begin{center}\footnotesize{\begin{tabular}{c|c} On $\MM$& On $\MMM$\\\hline
$\begin{array}{rll}
    \xx'^{0} = &\xx^{0} + \log J_{f}; &\\
    \xx'^{i} = &f^{i}(\xx^{1},\dots,\xx^{n}), &\mbox{for $i = 1,\dots,n$},\end{array}
$& $\begin{array}{rll}
    t' =&t J_{f}; &\\
    \xxx'^{i} =&f^{i}(\xxx^{1},\dots,\xxx^{n}) & \mbox{for $i = 1,\dots,n$},\end{array}
$\\
\end{tabular}}\end{center}
where $J_{f}$ denotes the modulus of the Jacobian of $f$ considered
as a function of coordinates.  We are immediately led to the
following observation.
\begin{theorem}\label{main}
\begin{enumerate}\item The bundles $\ppi:\MM\to M$ and $\pppi:\MMM\to M$ are isomorphic, that is, there is a diffeomorphism $F:\MMM\to \MM$ such that
$\ppi\circ F = \pppi$.
    \item The weight operator $\w$ is mapped to $\frac{\partial}{\partial (\xx^{0})}$ under $F$.
    \item A projective equivalence class $\Pi_{ij}^{k}$ on $M$ canonically induces a linear connection on
    $\MMM$, in particular allowing us to consider covariant derivatives of densities along vector
    fields on $M$.

\end{enumerate}
\end{theorem}
\begin{proof}
\begin{enumerate}
\item In coordinates, simply define
\[\xx^{i}(F(x)) = F^{i}(t, \xxx^{1},\dots,\xxx^{n}) = \left\{
\begin{array}{ll}
    \log t & \mbox{for $i = 0$, and $t>0$};\\
    \xxx^{i} & \mbox{for $i = 1,\dots,n$}.\end{array} \right.
\]
\item Written in terms of generating functions $\w$ takes the form
of a logarithmic derivative $t\frac{\partial}{\partial t}$, which is
mapped to $\frac{\partial}{\partial (\xx^{0})}$.
\item A projective equivalence class $\Pi^{k}_{ij}$ canonically defines, via Thomas' construction, a
linear connection $\nnabla$ on $\MM$ whose coefficients are given by
(\ref{t1})-(\ref{t3}).  Now define a linear connection on $\MMM$ by
pulling $\nnabla$ back along $F$.

\end{enumerate}
\end{proof}
\section{Operators, Brackets and future developments}

Let $M$ be a manifold endowed with a tensor field $S^{ij}$. This
defines a map $S^{\sharp}:T^{*}M\to TM$ given in local coordinates
by $\omega_{i}dx^{i}\mapsto \omega_{j}S^{ji}\partial_{i}$.

\begin{definition}Let $E\to M$ be a vector bundle over a manifold $M$ endowed with a
tensor field $S^{ij}$.  An \emph{upper connection} on $E$ over
$S^{ij}$ is an $\RR$-bilinear map $\nabla:\Omega^{1}(M)\times
\Gamma(E) \to \Gamma(E)$ satisfying
\begin{equation*}
\nabla^{f\omega}\sigma = f\nabla^{\omega}\sigma \mbox{  and  }
\nabla^{\omega}f\sigma = f\nabla^{\omega}\sigma + (S^{\sharp}\omega)(f)\sigma \mbox{ for }f\in C^{\infty}(M).
\end{equation*}
If $S^{ij}$ is an invertible matrix, upper and ordinary connections
correspond by raising and lowering indices.
\end{definition}

In \cite{oloii}, symmetric biderivations on the algebra of densities
(or more briefly \emph{brackets}) were considered. In a system of
local coordinates $x^{1},\dots x^{n}$, a homogeneous bracket
$\{\,\cdot\, ,\,\cdot\,\}$ is uniquely defined by a triple of
quantities $S^{ij}, \gamma^{i}, \theta$ via
\begin{equation*}
\{x^{i},x^{j}\} = S^{ij}(Dx)^{\lambda},\hspace{1cm} \{x^{i},Dx\} =
\gamma^{i}(Dx)^{\lambda + 1},\hspace{1cm} \{Dx,Dx\}
=\theta(Dx)^{\lambda + 2};
\end{equation*}
where $\lambda\in \RR$ is the \emph{weight} of the bracket.  Here
$S^{ij} = S^{ji}$ is symmetric.  For a bracket of weight zero,
$S^{ij}$ is a tensor and $\gamma^{i}$ defines an upper connection
over $S^{ij}$ in the bundle of volume forms on $M$.
\begin{theorem}Let $M$ be a manifold endowed with a projective equivalence class (Thomas projective connection $\Pi^{i}_{jk}$) and a tensor field
$S^{ij}$.  Then:
\begin{enumerate}\item These data define an upper connection over $S^{ij}$ in the bundle of volume forms given by the
coefficients
\begin{equation}\gamma^{i} = \frac{n+1}{n+3}(\partial_{j}S^{ij} + S^{jk}\Pi_{jk}^{i}).\label{upcon}\end{equation}
\item The following expression defines an invariant second order differential operator acting on functions on $M$
\begin{equation}
\Delta = S^{ij}\partial_{i}\partial_{j} +
\left(\frac{2}{n+3}\partial_{j}S^{ij} -
\left(\frac{n+1}{n+3}\right)S^{jk}\Pi^{i}_{jk}\right)\partial_{i}\,
.\label{op}
\end{equation}
\end{enumerate}
\end{theorem}
\begin{remark}
In general, a linear connection $\Gamma^{k}_{ij}$ gives rise to a
connection on the bundle of volume forms by taking the trace
$\gamma_{i} := \Gamma^{k}_{ki}$.  Suppose that $S^{ij}$ is
non-degenerate, i.e., $S_{ij} = g_{ij}$ defines a metric. An upper
connection may be obtained in two ways:
\begin{itemize}
\item Take the coefficients $\Gamma^{k}_{ij}$ of the Levi-Civita
connection associated with $g_{ij}$.  Define $\gamma^{i} :=
g^{ij}\Gamma^{k}_{kj}$.
\item Take the class $\Pi^{k}_{ij}$ associated with the Levi-Civita connection coefficients $\Gamma^{k}_{ij}$ and
use (\ref{upcon}).
\end{itemize}
These two constructions give exactly the same coefficients
$\gamma^{i}$.
\end{remark}
Given a Thomas projective connection on $M$, the principal symbol of
a second order differential operator (of arbitrary weight) on the
algebra of densities $\V(M)$ may similarly be extended to an
invariant operator on $\V(M)$. From a projective equivalence class
$\Pi^{k}_{ij}$ on $M$, Theorem \ref{main} gives a linear connection
$\gggamma^{k}_{ij}$ on $\MMM$. Applying formula (\ref{op}) to the
associated Thomas projective connection on $\MMM$ gives the result.

In the presence of a Thomas projective connection on $M$, this is a
method of constructing, from a bracket on $\V(M)$, a canonical
differential operator on $\V(M)$ generating the bracket. A similar
construction using a natural inner product was one of the main
results of \cite{oloii}. Comparing operators from (\ref{op}) and
\cite{oloii} gives expressions for $\gamma^{i}$ and $\theta$ in
terms of $S^{ij}$ and $\Pi^{k}_{ij}$. Given a manifold, we have
therefore a map from the space of symbols $S^{ij}$ to the space of
second order differential operators on the algebra of densities, or
using Ovsienko's terminology (see \cite{lecomte-1998}), a
projectively equivariant quantisation (in the non-flat case). Indeed
it would no doubt prove fruitful to investigate this relation
further.

\textbf{Acknowledgments:} The author is greatly indebted to H.
Khudaverdian and Th. Voronov for their invaluable guidance and many
discussions and most grateful to A. Odzijewicz for providing such a
wonderful forum in Bia{\l}owie\.{z}a in which to present the results
of this note.
\bibliography{bialpap}
\bibliographystyle{aipproc}
\end{document}